\documentclass[11pt]{amsart}
\usepackage{amssymb,latexsym}

\topskip  \numberwithin{equation}{section}
\newtheorem{theorem}{Theorem}[section]

\newtheorem{corollary}[theorem]{Corollary}

\newtheorem{lemma}[theorem]{Lemma}

\begin{document}

\pagenumbering{arabic}
\def\NN{\mathbb{N}}
\def\rises{\mbox{rises}}
\def\levels{\mbox{levels}}
\def\drops{\mbox{drops}}
\def\parts{\mbox{parts}}
\def\sof{$\Box$}
\def\ds{\displaystyle}
\def\Aodd{A=\{m\,|\,m=2k+1, k \ge 0\}}

\title{Counting rises, levels, and drops in compositions}
\maketitle
\begin{center}
Silvia Heubach and Toufik Mansour
\end{center}

\begin{center}
{Department of Mathematics, California State University Los Angeles\\
Los Angeles, CA 90032, USA}

   {\tt sheubac@calstatela.edu}\\[10pt]

{Department of Mathematics, Haifa University,
31905 Haifa, Israel

        {\tt toufik@math.haifa.ac.il} }
\end{center}
%
\section*{Abstract}
A composition of $n\in\NN$ is an ordered collection of one
or more positive integers whose sum is $n$. The number of summands
is called the number of  parts of the composition. A
palindromic composition of $n$ is a composition of $n$ in which
the summands are the same in the given or in reverse order. In
this paper we study the generating function for the number of
compositions (respectively palindromic compositions) of $n$ with $m$ parts in a
given set $A\subseteq\NN$ with respect to the number of rises,
levels, and drops. As a consequence, we derive all the
previously known results for this kind of problem, as well as
many new results.

\noindent {\bf AMS Classification Number}:  05A05, 05A15

\noindent {\bf Key words}: Compositions, palindromic compositions,
Carlitz compositions, partitions, generating functions.

\section{Introduction}

A {\em composition} $\sigma=\sigma_1\sigma_2\ldots\sigma_m$ of
$n\in\NN$ is an ordered collection of one or more positive
integers whose sum is $n$. The number of summands, namely $m$, is
called the number of {\em parts} of the composition. A {\em
palindromic composition} of $n\in\NN$ is a composition for which
$\sigma_1\sigma_2\ldots\sigma_m=
\sigma_m\sigma_{m-1}\ldots\sigma_1$. A {\em Carlitz composition}
is a composition of $n\in\NN$ in which no two consecutive parts
are the same. We will derive the generating functions for the
number of compositions, number of parts, and number of {\em rises}
(a summand followed by a larger summand), {\em levels} (a summand
followed by itself), and {\em drops} (a summand followed by a
smaller summand) in all compositions of $n$ whose parts are in a
given set $A$. This unified framework generalizes earlier work by
several authors.

Alladi and Hoggatt~\cite{AH} considered $A=\{1,2\}$, and derived
generating functions for the number of compositions, number of
parts, and number of rises, levels and drops in compositions and
palindromic compositions of $n$. Chinn and Heubach~\cite{CH2}
generalized to $A=\{1,k\}$ and derived all the respective
generating functions. Chinn, Grimaldi and Heubach~\cite{CGH}
considered the case $A=\NN$, and derived generating functions for
all quantities of interest. Grimaldi~\cite{G} studied
$A=\{m|m=2k+1,k\ge 0\}$, and derived generating functions for the
number of such compositions, as well as the number of parts, but
not for the number of rises, levels and drops. In addition, he
studied compositions without the summand 1~\cite{G2}, which was
generalized by Chinn and Heubach~\cite{CH}, who looked at
compositions without the summand $k$, i.e. $A=\NN-\{k\}$. In both
cases, the authors only derived generating functions for the total
number of compositions and the number of parts, but not for the
number of rises, levels and drops. Finally, Hoggatt and
Bricknell~\cite{HB} looked at compositions with parts in a general
set $A$, and gave generating functions for the number of
compositions and the number of parts. This work was generalized by
Heubach and Mansour~\cite{HM}, which also considered Carlitz
compositions and gave additional generating functions for the
number of compositions with a given number of parts in a set $B
\subseteq A$.

We will present a unified framework which allows us to derive
previous results by choosing a specific set $A$, as well as new
results. We will therefore study the specific sets $A=\NN$,
$A=\{1,2\}$, $A=\{1,k\}$, $A=\NN-\{k\}$, and
$A=\{m\,|\,m=2k+1,k\ge 0\}$. In the case of Carlitz compositions,
we will restrict ourselves to the sets $A=\{1,2\}$, $A=\{1,k\}$
and  $A=\{a,b\}$. The main result and its proof will be stated in
Section 2, and in Section~\ref{seccar} we present several
applications on the set of compositions (see
Subsection~\ref{sec31}), palindromic compositions (see
Subsection~\ref{sec32}), Carlitz compositions (see
Subsection~\ref{seccar1}), Carlitz palindromic compositions (see
Subsection~\ref{seccarp}), and partitions (see
Subsection~\ref{secpar}) of $n$ with $m$ parts in $A$,
respectively. As a consequence, we derive all the previously known
results for this kind of problem, as well as many new results.

\section{Main Result}
Let $\NN$ be the set of all positive integers, and let $A$ be any
ordered (finite or infinite) set of positive integers, say
$A=\{a_1,a_2,\ldots,a_k\}$, where $a_1<a_2<a_3<\cdots<a_k$, with
the obvious modifications in the case $|A|=\infty$. In the
theorems and proofs, we will treat the two cases together if
possible, and will note if the case $|A|=\infty$ requires
additional steps. For ease of notation, ``ordered set'' will
always refer to a set whose elements are listed in increasing
order.

For any ordered set $A=\{a_1,a_2,\ldots,a_k\}\subseteq\NN$, we
denote the set of all compositions (respectively palindromic
compositions) of $n$ with parts in $A$ by $C_n^A$ (respectively
$P_n^A$). For any composition $\sigma$, we denote the number of
parts, rises, levels, and drops by $\parts(\sigma)$,
$\rises(\sigma)$, $\levels(\sigma)$, and $\drops(\sigma)$,
respectively. We denote the generating function for the number of
compositions (respectively palindromic compositions) of $n$ with
$\parts(\sigma)$ parts in a set $A$ such that there are
$\rises(\sigma)$ rises, $\levels(\sigma)$ levels, and
$\drops(\sigma)$ drops by $C_A(x;y;r,\ell,d)$ (respectively
$P_A(x;y;r,\ell,d)$), that is,
$$C_A(x;y;r,\ell,d)=\sum_{n\geq0}\sum_{\sigma\in
C_n^A}x^ny^{\parts(\sigma)}r^{\rises(\sigma)}\ell^{\levels(\sigma)}d^{\drops(\sigma)}$$
and
$$P_A(x;y;r,\ell,d)=\sum_{n\geq0}\sum_{\sigma\in
P_n^A}x^ny^{\parts(\sigma)}r^{\rises(\sigma)}\ell^{\levels(\sigma)}d^{\drops(\sigma)}.$$

The main result of this paper gives explicit expressions for the generating functions $C_A(x;y;r,\ell,d)$
and $P_A(x;y;r,\ell,d)$.

\begin{theorem}\label{mth}
Let $A=\{a_1,\ldots,a_k\}$ be any ordered subset of $\NN$.

{\rm(i)} The generating function $C_A(x;y;r,\ell,d)$ is given by
$$\dfrac{1 +(1-d)\displaystyle{\sum}_{j=1}^k\left( \frac{\displaystyle x^{a_j}y}{\displaystyle1-x^{a_j}y(\ell-d)}{\displaystyle \prod_{i=1}^{j-1}}\frac{\displaystyle 1-x^{a_i}y(\ell-r)}{\displaystyle 1-x^{a_i}y(\ell-d)}\right)}
{1-d\displaystyle \sum_{j=1}^k\left(\frac{\displaystyle x^{a_j}y}{\displaystyle 1-x^{a_j}y(\ell-d)}\displaystyle \prod_{i=1}^{j-1}\frac{\displaystyle 1-x^{a_i}y(\ell-r)}{\displaystyle 1-x^{a_i}y(\ell-d)}\right)}.$$

{\rm(ii)} The generating function $P_A(x;y;r,\ell,d)$ is given by
$$\dfrac{1+\displaystyle \sum_{i=1}^k\frac{\displaystyle x^{a_i}y+x^{2a_i}y^2(\ell-d \, r)}{\displaystyle 1-x^{2a_i}y^2(\ell^2-d \, r)}}
{1-\displaystyle \sum_{i=1}^k\frac{\displaystyle x^{2a_i}y^2d \, r}{\displaystyle 1-x^{2a_i}y^2(\ell^2-d \, r)}}.$$
\end{theorem}

\subsection{Proof of Theorem~\ref{mth}(i)}
Our present aim is to find $C_A(x;y;r,\ell,d)$ explicitly, thus we
need the following definition. For all $e\geq1$, we define
$$C_A(s_1s_2\ldots s_e|x;y;r,\ell,d)=\sum_{n\geq0}\sum_{\sigma}x^ny^{\parts(\sigma)}r^{\rises(\sigma)}\ell^{\levels(\sigma)}d^{\drops(\sigma)},$$
where the sum on the right side of the equation is over all the
composition $\sigma\in C_n^A$ such that $\sigma_j=s_j$ for all
$j=1,2,\ldots,e$, i.e., the composition $\sigma$ starts with $s_1s_2\ldots s_e$.

Now, let us introduce two relations (Equation (\ref{eqaa}) and
Lemma~\ref{lemaa}) between the generating functions
$C_{A}(x;y;r,\ell,d)$ and $C_A(a_i|x;y;r,\ell,d)$.
The first relation is given by
\begin{equation}\label{eqaa}
C_A(x;y;r,\ell,d)=1+\sum_{i=1}^k C_A(a_i|x;y;r,\ell,d),
\end{equation}
which follows immediately from the definitions (note that the
summand 1 covers the case $n = 0$). The second relation is given
by the following lemma, and stems from a recursive creation of the
compositions of $n$.

\begin{lemma}\label{lemaa}
Let $A=\{a_1,\ldots,a_k\}$ be any ordered subset  of $\NN$.
For all $i=1,2,\ldots,k$, the generating function
$C_A(a_i|x;y;r,\ell,d)$ is given by
$$x^{a_i}y\left(1+d\sum_{j=1}^{i-1}C_A(a_j|x;y;r,\ell,d)+\ell
C_A(a_i|x;y;r,\ell,d)+r\sum_{j=i+1}^{k}C_A(a_j|x;y;r,\ell,d)\right).$$
\end{lemma}
\begin{proof}
The compositions of $n$ starting with $a_i$ with at least two
parts can be created recursively by prepending $a_i$ to a
composition of $n-a_i$ which starts with $a_j$ for some $j$. This
either creates a rise (if $i<j$), a level (if $i=j$), or a drop
(if $i > j$), and in each case, results in one more part. Thus,
$$C_A(a_ia_j|x;y;r,\ell,d)=\left\{
\begin{array}{ll}
rx^{a_i}yC_A(a_j|x;y;r,\ell,d),& i<j\\
\ell x^{a_i}yC_A(a_j|x;y;r,\ell,d),& i=j\\
dx^{a_i}yC_A(a_j|x;y;r,\ell,d),& i>j\\
\end{array}.
\right.$$ Summing over $j$ and accounting for the single composition with exactly one part,
namely $a_i$, gives the stated result.
\end{proof}

We are now ready to prove Theorem~\ref{mth}(i). Lemma~\ref{lemaa}
together with Equation~(\ref{eqaa}) results in a system of $k+1$
equations in $k+1$ variables, where we define
$t_{0}=C_A(x;y;r,\ell,d)$, $t_i=C_A(a_i|x;y;r,\ell,d)$ and
$b_i=x^{a_i}y$, for $i=1,2,\ldots,k$:
\begin{equation}\label{eqab}
{\small \left\{\begin{array}{l}
t_0-t_1-t_2-t_3\cdots-t_{k-1}-t_{k}=1\\
(1-b_1\ell)t_1-b_1rt_2-b_1rt_3\cdots-b_1rt_{k-1}-b_1rt_k=b1\\
-b_2dt_1+(1-b_2\ell)t_2-b_2rt_3\cdots-b_2rt_{k-1}-b_2rt_k=b2\\
-b_3dt_1-b_3dt_2+(1-b_3\ell)t_3\cdots-b_3rt_{k-1}-b_3rt_k=b3\\
\vdots\\
-b_{k-1}dt_1-b_{k-1}dt_2-b_{k-1}t_3\cdots+(1-b_{k-1}\ell)t_{k-1}-b_{k-1}rt_k=b_{k-1}\\
-b_{k}dt_1-b_{k}dt_2-b_{k}dt_3\cdots-b_kdt_{k-1}+(1-b_{k}\ell)t_k=b_{k}
\end{array}\right..}
\end{equation}

Let $M_k$ be the $(k+1)\times(k+1)$ matrix of the system of equations (\ref{eqab}), i.e.,
{\small$$M_k=\left(\begin{array}{ccccccc}
1&-1&-1&\cdots&-1&-1&-1\\
0&1-b_1\ell&-b_1r&-b_1r&\cdots&-b_1r&-b_1r\\
0&-b_2d&1-b_2\ell&-b_2r&\cdots&-b_2r&-b_2r\\
0&-b_3d&-b_3d&1-b_3\ell&\cdots&-b_3r&-b_3r\\
\vdots&&&\vdots&&&\vdots\\
0&-b_{k-1}d&-b_{k-1}d&-b_{k-1}&\cdots&1-b_{k-1}\ell&-b_{k-1}r\\
0&-b_{k}d&-b_{k}d&-b_{k}d&\cdots&-b_kd&1-b_{k}\ell
\end{array}\right).$$} We also define the $(k+1)\times(k+1)$ matrix $N_k$ which results
from replacing the first column in $M_k$ by the vector of the
right-hand side of (\ref{eqab}), i.e.,
{\small$$N_k=\left(\begin{array}{ccccccc}
1&-1&-1&\cdots&-1&-1&-1\\
b_1&1-b_1\ell&-b_1r&-b_1r&\cdots&-b_1r&-b_1r\\
b_2&-b_2d&1-b_2\ell&-b_2r&\cdots&-b_2r&-b_2r\\
b_3&-b_3d&-b_3d&1-b_3\ell&\cdots&-b_3r&-b_3r\\
\vdots&&&\vdots&&&\vdots\\
b_{k-1}&-b_{k-1}d&-b_{k-1}d&-b_{k-1}&\cdots&1-b_{k-1}\ell&-b_{k-1}r\\
b_k&-b_{k}d&-b_{k}d&-b_{k}d&\cdots&-b_kd&1-b_{k}\ell
\end{array}\right).$$}
Then, by Cramer's Rule,
$t_{0}=C_A(x;y;r,\ell,d)=\frac{\det(N_k)}{\det(M_k)}$. We now
derive formulas for these two determinants. Expanding down the
first column of $M_k$, we get that
{\small$$\det(M_k)=\left|\begin{array}{llllll}
1-b_1\ell&-b_1r&-b_1r&\cdots&-b_1r&-b_1r\\
-b_2d&1-b_2\ell&-b_2r&\cdots&-b_2r&-b_2r\\
-b_3d&-b_3d&1-b_3\ell&\cdots&-b_3r&-b_3r\\
\vdots&&&\vdots&&\vdots\\
-b_{k-1}d&-b_{k-1}d&-b_{k-1}&\cdots&1-b_{k-1}\ell&-b_{k-1}r\\
-b_{k}d&-b_{k}d&-b_{k}d&\cdots&-b_kd&1-b_{k}\ell
\end{array}\right|.$$}
Subtracting the $(k-1)^{\mbox{\footnotesize st}}$ column from
$k^{\mbox{\footnotesize th}}$ column of the above matrix, then
expanding down the resulting column gives that
\begin{equation}\label{detM}
\det(M_k)=
(1-b_k(\ell-d))\det(M_{k-1})-b_kd(1-b_{k-1}(\ell-r))\det(E(b_1,b_2,...,b_{k-2})),
\end{equation}
where
$$E(b_1,b_2,...,b_{k-2})=
\left(\begin{array}{ccccc}
1-b_1\ell&-b_1r&-b_1r&\cdots&-b_1r\\
-b_2d&1-b_2\ell&-b_2r&\cdots&-b_2r\\
-b_3d&-b_3d&1-b_3\ell&\cdots&-b_3r\\
\vdots&&\vdots&&\vdots\\
-b_{k-2}d&-b_{k-2}d&-b_{k-2}d&\cdots&-b_{k-2}r\\
1&1&1&\cdots&1
\end{array}\right).$$
Adding $(b_1r)$ times the last row to the first row in the matrix
$E(b_1,b_2,...,b_{k-2})$, then expanding across the resulting
first row gives
$$\det(E(b_1,b_2,...,b_{k-2}))=(1-b_1(\ell-r))\det(E(b_2,...,b_{k-2})),$$
and, since $\det(E(b_{k-2}))=(1-b_{k-2}(\ell-r))$,
\begin{equation}\label{detE}
\det(E(b_1,b_2,...,b_{k-2}))=\prod_{j=1}^{k-2}(1-b_j(\ell-r)).\end{equation}
Equations~(\ref{detM}) and (\ref{detE}) result in
$$\det(M_k)=
(1-b_k(\ell-d))\det(M_{k-1})-b_kd\prod_{j=1}^{k-1}(1-b_j(\ell-r))
.$$
Thus, if we define $\det(M_0)=1$ and use the fact that
$\det(M_1)=1-b_1\ell=1-b_1(\ell-d)-b_1d$, then we can show by
induction on $k$ that for all $k\geq1$,
\begin{equation}\label{eqac}
\det(M_k)=\prod_{j=1}^k(1-b_j(\ell-d))-d\sum_{j=1}^k
b_j\prod_{i=1}^{j-1}(1-b_i(\ell-r))\prod_{i=j+1}^{k}(1-b_j(\ell-d)).
\end{equation}
Similarly, by subtracting $(b_kd)$ times the last row from the
$k^{\mbox{\footnotesize th}}$ row in the matrix $N_k$ and then
expanding across the resulting $k^{\mbox{\footnotesize th}}$ row
we get
\begin{equation}\label{detN}
\det(N_k)=
(1-b_k(\ell-d))\det(N_{k-1})+b_k(1-d)\det(D(b_1,b_2,...,b_{k-1})),
\end{equation}
where $D(b_1,b_2,...,b_{k-1})$ agrees with
$E(b_1,b_2,...,b_{k-1})$ except for the signs of the last row.
Thus,
$\det(D(b_1,b_2,...,b_{k-1}))=-\det(E(b_1,b_2,...,b_{k-1}))$,
which yields
$$\det(N_k)=
(1-b_k(\ell-d))\det(N_{k-1})-b_k(1-d)\prod_{j=1}^{k-1}(1-b_j(\ell-r)).$$
With $\det(N_0)=1$ and
$\det(N_1)=1-b_1\ell+b_1=1-b_1(\ell-d)+(1-d)b_1$, we can show by
induction on $k$ that for all $k\geq1$,
\begin{equation} \label{eqac2}
\det(N_k)=\prod_{j=1}^k(1-b_j(\ell-d))+(1-d)\sum_{j=1}^k
b_j\prod_{i=1}^{j-1}(1-b_i(\ell-r))\prod_{i=j+1}^{k}(1-b_j(\ell-d)).
\end{equation}
Substituting Equations~(\ref{eqac}) and (\ref{eqac2}) and
$b_i=x^{a_i}y$ into $\frac{\det(N_k)}{\det(M_k)}$ completes the
proof of Theorem~\ref{mth}(i). Note that if $|A|=\infty$, then the
result follows by taking limits as $k \rightarrow \infty$. \sof
\subsection{Proof of Theorem~\ref{mth}(ii)}
As in the proof of part (i), we need to find an explicit
expression for $P_A(x;y;r,\ell,d)$, thus we define for all
$e\geq1$
$$P_A(s_1s_2\ldots s_e|x;y;r,\ell,d)=\sum_{n\geq0}\sum_{\sigma}x^ny^{\parts(\sigma)}r^{\rises(\sigma)}\ell^{\levels(\sigma)}d^{\drops(\sigma)},$$
where the sum on the right side of the equation is over all the
palindromic compositions $\sigma\in P_n^A$ such that
$\sigma_j=s_j$ for all $j=1,2,\ldots,e$.

As before, we get two relations (Equation (\ref{eqba}) and
Lemma~\ref{lemba}) between the generating functions
$P_{A}(x;y;r,\ell,d)$ and $P_A(a_i|x;y;r,\ell,d)$.
The first relation is given by
\begin{equation}\label{eqba}
P_A(x;y;r,\ell,d)=1+\sum_{i=1}^k P_A(a_i|x;y;r,\ell,d),
\end{equation}
which holds immediately from the definitions. The second relation
is given by the following lemma.

\begin{lemma}\label{lemba}
Let $A=\{a_1,\ldots,a_k\}$ be any ordered subset of $\NN$. For all
$i=1,2,\ldots,k$, the generating function $P_A(a_i|x;y;r,\ell,d)$
is given by
$$\frac{x^{a_i}y+x^{2a_i}y^2\ell}{1-x^{2a_i}y^2(\ell^2-d \, r)}
+\frac{x^{2a_i}y^2d \, r}{1-x^{2a_i}y^2(\ell^2-d \, r)}(P_A(x;y;r,\ell,d)-1).$$
\end{lemma}
\begin{proof}
First of all, for all $e,m\geq1$ we define
$$P_A(s_1s_2\ldots s_e;m|x;r,\ell,d)=\sum_{n\geq0}\sum_{\sigma}x^nr^{\rises(\sigma)}\ell^{\levels(\sigma)}d^{\drops(\sigma)},$$
where the sum on the right side of the equation is over all the
palindromic compositions $\sigma\in P_n^A$ with $m$ parts such
that $\sigma_j=s_j$ for all $j=1,2,\ldots,e$.

Now, by fixing $i$ and using the definitions we have that
$$P_A(a_i;0|x;r,\ell,d)=0,\ P_A(a_i;1|x;r,\ell,d)=x^{a_i},\
P_A(a_i;2|x;r,\ell,d)=x^{2a_i}\ell.$$ For $m\geq3$, we create the
palindromic compositions of $n$ that start and end with $a_i$ from
those of $n-2 a_i$ that start with $a_j$ by prepending and
appending $a_i$. This results in two additional parts, and in one
additional drop and rise when $i \neq j$, and two additional
levels when $i=j$. (Note that the symmetry of the palindromic
compositions, which distinguishes only the case $i=j$, allows us
to use a different proof technique, which does not work for
compositions.) Thus, for $m\geq3$,
$$\begin{array}{l}
P_A(a_i;m|x;r,\ell,d)=\sum_{j=1,\,j\neq i}^kP_A(a_ia_j;m|x;r,\ell,d)+P_A(a_ia_i;m|x;r,\ell,d)\\
\\[-4pt]
\qquad=x^{2a_i}d\, r\sum_{j=1,\,j\neq i}^kP_A(a_j;m-2|x;r,\ell,d)+x^{2a_i}\ell^2P_A(a_i;m-2|x;r,\ell,d)\\
\\[-4pt]
\qquad=x^{2a_i}d\, r\sum_{j=1}^k
P_A(a_j;m-2|x;r,\ell,d)+x^{2a_i}(\ell^2-d\, r)P_A(a_i;m-2|x;r,\ell,d).
\end{array}$$

Multiplying by $y^m$ and summing over all $m\geq0$,
we get that
$$\begin{array}{l}
P_A(a_i|x;y;r,\ell,d)=x^{a_i}y+x^{2a_i}y^2\ell+x^{2a_i}y^2d\,r\sum_{j=1}^k P_A(a_j|x;y;r,\ell,d)\\
\qquad\qquad\,\,\,\,\,\qquad\qquad\qquad+x^{2a_i}y^2(\ell^2-d\,r)P_A(a_i|x;y;r,\ell,d),
\end{array}$$ or, equivalently,
$$P_A(a_i|x;y;r,\ell,d)=\frac{x^{a_i}y+x^{2a_i}y^2\ell}{1-x^{2a_i}y^2(\ell^2-d\,r)}
+\frac{x^{2a_i}y^2d\,r}{1-x^{2a_i}y^2(\ell^2-d\,r)}\sum_{j=1}^k
P_A(a_j|x;y;r,\ell,d),$$ from which Lemma~\ref{lemba} follows by
using Equation~(\ref{eqba}).
\end{proof}

Now we are ready to give the proof of Theorem~\ref{mth}(ii).
Applying Lemma~\ref{lemba} for all $i=1,2,\ldots,k$ together with
using Equation~(\ref{eqba}), we get that the generating function
$P_A(x;y;r,\ell,d)$ is given by
$$1+\sum_{i=1}^k\frac{x^{a_i}y+x^{2a_i}y^2(\ell-d\,r)}{1-x^{2a_i}y^2(\ell^2-d\,r)}
+\sum_{i=1}^k\frac{x^{2a_i}y^2d\,r}{1-x^{2a_i}y^2(\ell^2-d\,r)}P_A(x;y;r,\ell,d).$$
Equivalently,
$$P_A(x;y;r,\ell,d)=\dfrac{1+\sum\limits_{i=1}^k\dfrac{x^{a_i}y+x^{2a_i}y^2(\ell-dr)}{1-x^{2a_i}y^2(\ell^2-dr)}}
{1-\sum\limits_{i=1}^k\dfrac{x^{2a_i}y^2dr}{1-x^{2a_i}y^2(\ell^2-dr)}},$$
as claimed. \sof
\section{Applications}\label{seccar}
In the following subsections we give several applications for both parts of Theorem~\ref{mth}.

\subsection{Compositions with parts in $A$}\label{sec31}
In this subsection we study the number of compositions of $n$ as well as the number
of rises, levels, and drops in the compositions of $n$ with parts in $A$.
Applying Theorem~\ref{mth}(i)
for $r=1$, $\ell=1$, and $d=1$, we get that the generating function
for the number of compositions of $n$ with $m$ parts in $A$ is
given by
\begin{equation} \label{comp-n-m}
           \dfrac{1}{1-y\sum_{j=1}^kx^{a_j}}.
\end{equation}
Therefore, the generating function for the number of compositions
of $n$ with $m$ parts in $\NN$ is given by
$$\sum_{n\geq0}\sum_{\sigma\in
C_n^{\NN}}x^ny^{\parts(\sigma)}=\frac{1}{1-y\sum_{j=1}^{\infty}x^j}=\frac{1}{1-\frac{y\,x}{1-x}}=\sum_{m\geq0}\frac{x^m}{(1-x)^m}y^m.$$

Furthermore, setting $y=1$ in Equation~(\ref{comp-n-m}) gives the
generating function for the number of compositions of $n$ with
parts in $A$ (see~\cite{HM}, Theorem 2.4):
$$ \dfrac{1}{1-\sum_{j=1}^kx^{a_j}}.$$
In particular, for $A=\NN$, the generating function for the number
of compositions of $n$ with parts in $\NN$ is given by (see \cite{CGH},
Theorem 6)
$$\frac{1-x}{1-2x}.$$

Additional examples for specific choices of $A$ are given in~\cite{HM}.

\subsubsection{Number of rises and drops}
Note that the number of rises always equals the number of drops in
all compositions of $n$: for each non-palindromic composition
there exists a composition in reverse order, thus the rises match
the drops, and for palindromic compositions, symmetry matches up
rises and drops within the composition. Thus, we will derive
results only for rises, and the results for drops follow by
interchanging the roles of $r$ and $d$ in the proofs.

Setting $\ell=1$ and $d=1$ in Theorem~\ref{mth}(i) gives
\begin{equation}\label{eqrr}
C_A(x;y;r,1,1)=\dfrac{1}{1-\sum_{j=1}^k\left(x^{a_j}y\prod_{i=1}^{j-1}(1-x^{a_i}y(1-r))\right)}.
\end{equation}
Using Equation~(\ref{eqrr}) together with the fact that for $f_i(r) \ne 0$
\begin{equation}\label{eqdiff}
\frac{\partial}{\partial r}\prod_{i=1}^m f_i(r)=\left(\prod_{i=1}^m f_i(r)\right)\sum_{i=1}^m\frac{\frac{\partial}{\partial r}f_i(r)}{f_i(r)},
\end{equation}
we get that
\begin{equation}\label{eqnrr}
\left.\frac{\partial}{\partial r}
C_A(x;y;r,1,1)\right|_{r=1}=\frac{y^2\sum_{k\geq j>i\geq 1}
x^{a_i+a_j}}{\left(1-y\sum_{j=1}^kx^{a_j}\right)^2}.
\end{equation}

Hence, expressing this function as a power series about $y=0$, we
get the following result.

\begin{corollary}\label{co2}
Let $A=\{a_1,\ldots,a_k\}$ be any ordered subset of $\NN$. Then
$$\sum_{n\geq0}\sum_{\sigma\in
C_n^A}\rises(\sigma)x^ny^{\parts(\sigma)}=\left(\sum_{k\geq j>i\geq 1}
x^{a_i+a_j}\right)\sum_{m\geq0}(m+1)\left(\sum_{j=1}^kx^{a_j}\right)^my^{m+2}$$
and
$$\sum_{n\geq0}\sum_{\sigma\in
C_n^A}\drops(\sigma)x^ny^{\parts(\sigma)}=\left(\sum_{k\geq j>i\geq 1}
x^{a_i+a_j}\right)\sum_{m\geq0}(m+1)\left(\sum_{j=1}^kx^{a_j}\right)^my^{m+2}.$$
\end{corollary}

For example, letting $A=\NN$ and looking at the coefficient of
$y^m$ in Corollary~\ref{co2} we get that the generating function
for the number of rises (drops) in the compositions of $n$ with a
fixed number of parts, $m\geq2$,  in $\NN$ is given by
\begin{eqnarray*}
\lefteqn{\sum_{j>i\geq1}x^{i+j}(m-1)(\sum_{j\geq1}x^j)^{m-2}} \hspace{.5in}\\
&=&\sum_{i\geq1}\sum_{j\geq i+1}
x^{i+j}(m-1)\left(\frac{x}{1-x}\right)^{m-2}=
\sum_{i\ge1}x^i\sum_{j\ge1}x^{2j}\frac{(m-1)
x^{m-2}}{(1-x)^{m-2}}\\
&=&\frac{x^3}{(1-x)(1-x^2)}\cdot \frac{(m-1)
x^{m-2}}{(1-x)^{m-2}}=\frac{(m-1)x^{m+1}}{(1+x)(1-x)^{m}}.
\end{eqnarray*}

Furthermore, setting $y=1$ and $A=\NN$ in Corollary~\ref{co2}
allows us to compute the generating function for the number of rises
(drops) in all compositions of $n$ with parts in $\NN$ (see \cite{CGH}, Theorem
6) in a similar way:
\begin{eqnarray*}\sum_{n\geq0}\sum_{\sigma\in
C_n^A}\rises(\sigma)x^n&=&\sum_{j>i\geq 1}
x^{i+j}\sum_{m\geq0}(m+1)\left(\frac{x}{1-x}\right)^m=\frac{x^3}{(1-x)(1-x^2)}\cdot \frac{1}{(1-\frac{x}{1-x})^2}\\
&=&\frac{x^3}{(1+x)(1-2x)^2}.
\end{eqnarray*}

In other words, as shown in ~\cite[Theorem~3]{CGH}, the number of
rises (drops) in the compositions of $n$ with parts in $\NN$ is
given by
$$\frac{1}{9}\left(2^{n-2}(3n-5)+(-1)^{n+1}\right) \qquad \mbox{for  } n \ge 3.$$

For $A=\{1,k\}$ and $y=1$, Corollary~\ref{co2} gives the
generating function for the number of rises (drops) in all
compositions of $n$ with parts in $\{1,k\}$ as (see~\cite{CH2},
Theorem 4)
$$\frac{x^{k+1}}{(1-x-x^k)^2}.$$

For $\Aodd$ and $y=1$, and using that $\sum_{0 \le i
<j}x^{(2i+1)+(2j+1)}=\sum_{i\ge 0}(x^2)^i \sum_{j \ge 1}(x^4)^j$,
Corollary~\ref{co2} yields a new result, namely that the
generating function for the number of rises (drops) in
compositions of $n$ with odd parts is given by
$$\frac{x^{k+1}}{(1-x-x^k)^2}.$$

For $A=\NN-\{k\}$, and defining
$g(x,y;k)=\sum_{n\geq0}\sum_{\sigma\in
C_n^A}\rises(\sigma)x^ny^{\parts(\sigma)}$, Corollary~\ref{co2}
gives
$$g(x,y;k)=\left(\frac{x^3}{(1-x)(1-x^2)}-\frac{x^{k+1}(1-x^{k-1})+x^{2k+1}}{1-x}\right)\sum_{m\geq0}(m+1)
\left(\frac{x}{1-x}-x^k\right)^my^{m+2}.$$
For $k=1$, i.e., $A=\NN-\{1\}$ we get that
$$g(x,y;1)
=\sum_{m\geq2}(m-1)\frac{x^{2m+1}}{(1+x)(1-x)^m}y^m.$$ Thus, the
generating function for the number of rises (drops) in the
compositions of $n$ with a fixed number of parts, $m\geq2$, in $A=\NN-\{1\}$ is
given by
$$(m-1)\frac{x^{2m+1}}{(1+x)(1-x)^m}=\sum_{n\geq0}x^{n+2m-1}(m-1)\sum_{j=0}^n(-1)^{n-j}\binom{j+m-1}{m-1}.$$

\subsubsection{Number of levels} Theorem~\ref{mth}(i) for $r=1$ and
$d=1$ gives
\begin{equation}\label{eqll}
C_A(x;y;1,\ell,1)=\dfrac{1}{1-\sum_{j=1}^k\frac{x^{a_j}y}{1-x^{a_j}y(\ell-1)}}.
\end{equation}
Therefore, using Equation (\ref{eqll}) we have that
\begin{equation}\label{eqnll}
\left.\frac{\partial}{\partial \ell}
C_A(x;y;1,\ell,1)\right|_{\ell=1}=\frac{y^2\sum_{j=1}^k
x^{2a_j}}{\left(1-y\sum_{j=1}^kx^{a_j}\right)^2}.
\end{equation}
Expressing the above function as a power series about $y=0$, we get the
following result.

\begin{corollary}\label{co3}
Let $A=\{a_1,\ldots,a_k\}$ be any ordered subset of $\NN$. Then
$$\sum_{n\geq0}\sum_{\sigma\in
C_n^A}\levels(\sigma)x^ny^{\parts(\sigma)}=\left(\sum_{j=1}^k
x^{2a_j}\right)\sum_{m\geq0}(m+1)\left(\sum_{j=1}^kx^{a_j}\right)^my^{m+2}.$$
\end{corollary}

Using computations similar to those for rises and drops, by
looking at the coefficient of $y^m$, we get from
Corollary~\ref{co3} that the generating function for the number of
levels in all compositions of $n$ with a fixed number of parts $m$ in
$\NN$ is given by
            $$\frac{(m-1)x^{m}}{(1+x)(1-x)^{m-1}}.$$
In addition, by setting $y=1$ and $A=\NN$ in
Corollary~\ref{co3} we obtain that the generating function for the
number of levels in the compositions of $n$ with parts in $\NN$ (see \cite{CGH}, Theorem 6) is given
by
$$\frac{x^2(1-x)}{(1+x)(1-2x)^2}.$$
Thus, as shown in ~\cite[Theorem~3]{CGH}, the number of levels in all compositions of $n$ with
parts in $\NN$  is given by
$$\frac{1}{9}\left(2^{n-2}(3n+1)+2(-1)^{n}\right)\qquad \mbox{for  } n \ge 1.$$

Applying Corollary~\ref{co3} for $A=\{1,2\}$ and $y=1$, we get the
result given in Theorem 1.1~\cite{AH} for the generating function
for the number of levels in all compositions with only 1's and
2's:
$$\frac{x^2+x^4}{(1-(x+x^2))^2},$$
and more generally, for $A=\{1,k\}$ and $y=1$, we get the result stated in Theorem 4~\cite{CH2}:
$$\frac{x^2+x^{2k}}{(1-(x^k+x^{2k}))^2}.$$

If we apply Corollary~\ref{co3} to $\Aodd$, then we get a new
result, namely that the generating function for the number of
levels in the compositions of $n$ with odd summands is given by
$$\frac{x^2(1-x^2)}{(1+x^2)(1-x-x^2)^2}.$$

Finally, we look at $A=\NN-\{k\}$ and define
$g(x,y;k)=\sum_{n\geq0}\sum_{\sigma\in
C_n^A}\levels(\sigma)x^ny^{\parts(\sigma)}$. Then
Corollary~\ref{co3} gives
$$g(x,y;k)=\left(\frac{x^2}{1-x^2}-x^{2k}\right)\sum_{m\geq0}(m+1)
\left(\frac{x}{1-x}-x^k\right)^my^{m+2}.$$
If we set $y=1$, then we get a new result, namely that the generating function for the number of levels in the compositions of $n$ without $k$ is given by
$$\frac{(1-x)x^2(1-x^{2(k-1)}+x^{2k})}{(1+x)(1-2x+x^k-x^{k+1})^2}.$$

\subsection{Palindromic compositions with parts in
$A$}\label{sec32} Applying Theorem~\ref{mth}(ii) for $r=1$,
$\ell=1$, and $d=1$ we get that the generating function for the
number of palindromic compositions of $n$ with $m$ parts in $A$ is
given by
$$\frac{1+y\sum_{i=1}^kx^{a_i}}{1-y^2\sum_{i=1}^k x^{2a_i}}.$$
Setting $y=1$ we get that the number of palindromic
compositions of $n$ with parts in $A$ is given by (see \cite{HM}, Theorem 3.2)
$$\frac{1+\sum_{i=1}^kx^{a_i}}{1-\sum_{i=1}^k x^{2a_i}}.$$
Using $A=\NN$ we get that the generating
function for the number of palindromic compositions of $n$ with
parts in $\NN$ is given by (see \cite{CGH}, Theorem 6)
$$\frac{1+x}{1-2x^2}.$$
Therefore, the number of palindromic compositions of $n$ with parts in
$\NN$ is given by  $2^{\lfloor n/2 \rfloor}$ (see \cite{CGH}, Theorem 1).

\subsubsection{{\bf Number of rises or drops}} As before,
the number of rises equals the number of drops.
Theorem~\ref{mth}(ii) for $\ell=1$ and $d=1$ gives
$$P_A(x;y;r,1,1)=\dfrac{1+\sum\limits_{i=1}^k\dfrac{x^{a_i}y+x^{2a_i}y^2(1-r)}{1-x^{2a_i}y^2(1-r)}}
{1-\sum\limits_{i=1}^k\dfrac{x^{2a_i}y^2r}{1-x^{2a_i}y^2(1-r)}}.$$
Therefore, by finding $\frac{\partial }{\partial r}
P_A(x;y;r,1,1)$ and setting $r=1$ we obtain the following
result.


\begin{corollary}\label{cp1}
Let $A=\{a_1,\ldots,a_k\}$ be any ordered subset of $\NN$. Then
the generating function  $g_A(x;y)=\sum_{n\geq0}\sum_{\sigma\in
P_n^A}\rises(\sigma)x^ny^{\parts(\sigma)}=\sum_{n\geq0}\sum_{\sigma\in
P_n^A}\drops(\sigma)x^ny^{\parts(\sigma)}$ is given by
$$\frac{y^2\left(1+y\sum_{i=1}^kx^{a_i}\right)\left(\sum_{i=1}^kx^{2a_i}(1-x^{2a_i}y^2)\right)
-y^2\left(1-y^2\sum_{i=1}^kx^{2a_i}\right)\sum_{i=1}^kx^{2a_i}(1+x^{a_i}y)}{\left(1-y^2\sum_{i=1}^k
x^{2a_i}\right)^2}.$$
\end{corollary}

For example, if $A=\NN$, then Corollary~\ref{cp1} gives that
\begin{equation} \label{palN}
g_{\NN}(x;y)=\frac{y^2\left(1+\frac{yx}{1-x}\right)\left(\frac{x^2}{1-x^2}-\frac{x^4y^2}{1-x^4})\right)
-y^2\left(1-\frac{y^2x^2}{1-x^2}\right)\left(\frac{x^2}{1-x^2}+\frac{x^3y}{1-x^3}\right)}{\left(1-\frac{y^2x^2}{1-x^2}\right)^2}.
\end{equation}
Thus, we can derive the generating function for the number of
rises (drops) in the compositions of $n$ with a given number of
parts, $m$, in $\NN$, by looking at the coefficient of $y^m$ in
$g_{\NN}(x;y)$. To do so, we expand the numerator of
$g_{\NN}(x;y)$ and collect terms according to powers of $y$:
$$\frac{x^4y^3}{(1-x)^2(1+x)} \cdot \left(\frac{2x+1}{(x^2+x+1)}+\frac{2x^2}{(x+1)(x^2+1)}y-\frac{
x^2}{(x^2+x+1)(x^2+1)}y^2\right).$$ Furthermore,
$$\frac{1}{\left(1-\frac{y^2x^2}{1-x^2}\right)^2}=\sum_{m \ge
0}\frac{(m+1)x^{2m}}{(1-x^2)^m}y^{2m},$$ so altogether,
\begin{eqnarray*}
g_{\NN}(x;y)&=&\sum_{m\geq0}\frac{(m+1)x^{2m+4}}{(1-x)^2(1+x)(1-x^2)^m}y^{2m+3}\cdot \\
&&\left(\frac{2x+1}{(x^2+x+1)}+\frac{2x^2}{(x+1)(x^2+1)}y-\frac{
x^2}{(x^2+x+1)(x^2+1)}y^2\right).
\end{eqnarray*}
We now have to distinguish between two cases, namely, $m$ odd and $m$
even. In the first case, only the summand with factor $y$ needs to
be taken into account, whereas in the second case, the summands
with factors $y^0$ and $y^2$ need to be considered. Thus, the
generating function for the number of rises (drops) in the
compositions of $n$ with a given number of parts, $m$, in $\NN$ is
given by
            $$\frac{(2m'-2)x^{2m'+2}}{(1+x^2)(1-x^2)^{m'}} \qquad \mbox{for  } m=2m',$$
and
            $$\frac{x^{2m'}(1-x)(1+(2m'-2)x+(2m'-3)x^2+(2m'-2)x^3)}{(1+x^2)(1+x+x^2)(1-x^2)^{m'}} \qquad \mbox{for  } m=2m'-1.$$
Furthermore, setting $y=1$ in Equation~(\ref{palN}) and simplifying
 yields that the generating function for
the number of rises (drops) in the compositions of $n$ with parts in $\NN$
(see~\cite{CGH}, Theorem~6) is given by
    $$g_{\NN}(x;1)=\frac{x^4(4x^4+4x^3+4x^2+3x+1)}{(1+x^2)(1+x+x^2)(1-2x^2)^2}.$$

We now apply Corollary~\ref{cp1} for $A=\{1,k\}$ and get that
$$g_{\{1,k\}}(x;y)=\frac{x^{k+1}y^3(x+x^k+2x^{k+1}y-y^2(x^3+x^{3k}-x^{k+2}-x^{2k+1}))}
{(1-y^2(x^2+x^{2k}))^2}.$$
In particular, when setting $y=1$ in
the above expression we get that the generating function for the
number of rises (drops) in the palindromic compositions of $n$
with any number of parts in $A=\{1,k\}$ is given by
(see~\cite{CH2}, Theorem 5)
$$g_{\{1,k\}}(x;1)=\frac{x^{k+1}(x-x^3+x^k-x^{3k}+2x^{k+1}+x^{k+2}+x^{2k+1})}
{(1-x^2-x^{2k})^2}.$$

If we let $\Aodd$ in Corollary~\ref{cp1}, then we get that the
generating function $g_{A}(x;y)$ is given by
$$\dfrac{y^2\left(1+\frac{x\,y}{1-x^2}\right)\left(\frac{x^2}{1-x^4}-\frac{y^2x^4}{1-x^8}\right)-y^2\left(1-\frac{x^2y^2}{1-x^4}\right)\left(\frac{x^2}{1-x^4}+\frac{x^3y}{1-x^6}\right)}{\left(1-\frac{x^2y^2}{1-x^4}\right)^2}.$$
Furthermore, if we let $y=1$ in the above expression, then we get
that the generating function for the number of rises (drops) in
the palindromic compositions of $n$ with any number of odd parts
is given by
$$g_A(x;1)=\frac{x^5(1+2x^2+2x^3+2x^4+2x^5+3x^6+2x^7+2x^8)}{(1+x^4)(1-x^2-x^4)^2(1+x^2+x^4)},$$
which extends the work of Grimaldi~\cite{G}.

Applying Corollary~\ref{cp1} to $A=\NN-\{k\}$ gives that
\begin{eqnarray*}g_{\NN-\{k\}}(x;y)&=&\frac{y^2\left(1+\frac{yx}{1-x}-yx^k\right)\left(\frac{x^2}{1-x^2}-x^{2k}-\frac{y^2x^4}{1-x^4}+y^2x^{4k}\right)}{\left(1-\frac{y^2x^2}{1-x^2}+y^2x^{2k}\right)^2}\\
&&-\,\frac{y^2\left(1-\frac{y^2x^2}{1-x^2}+y^2x^{2k}\right)\left(\frac{x^2}{1-x^2}-x^{2k}+\frac{yx^3}{1-x^3}-yx^{3k}\right)}
{\left(1-\frac{y^2x^2}{1-x^2}+y^2x^{2k}\right)^2}.
\end{eqnarray*} In
particular, when setting $y=1$ in the above expression we get that
the generating function for the number of rises (drops) in the
palindromic compositions of $n$ with any number of parts in
$A=\NN-\{k\}$ is given by
\begin{eqnarray*}
g_{\NN-\{k\}}(x;1)&=&\frac{x^4(1+3x+4x^2+4x^3+4x^4)+x^{2k+1}(x^4-1)(1+4x+5x^2+4x^3)}{(1+x^2)(1+x+x^2)(1-2x^2+x^{2k}-x^{2(k+1)})^2}\\&&+\frac{(x^2-1)(x^{k+2}+x^{3k}(1+x^2)(3x^2-2)+x^{4k}(1+x)(x-2))}{(1+x^2)(1-2x^2+x^{2k}-x^{2(k+1)})^2}.
\end{eqnarray*}
This extends the work of Chinn and Heubach~\cite{CH}. Likewise, we
can extend the work of Grimaldi~\cite{G2} by setting $k=1$ to get that
$$g_{\NN-\{1\}}(x;1)=\frac{(x^5+3x^4+5x^3+3x^2+3x+1)x^7}{(1-x^2-x^4)^2(1+x+x^2)(1+x^2)}.$$

\subsubsection{Number of levels} Theorem~\ref{mth}(ii) for $r=1$ and
$d=1$ gives
$$P_A(x;y;1,\ell,1)=\dfrac{1+\sum\limits_{i=1}^k\dfrac{x^{a_i}y+x^{2a_i}y^2(\ell-1)}{1-x^{2a_i}y^2(\ell^2-1)}}
{1-\sum\limits_{i=1}^k\dfrac{x^{2a_i}y^2}{1-x^{2a_i}y^2(\ell^2-1)}}.$$
Therefore, finding $\frac{\partial }{\partial \ell}
P_A(x;y;1,\ell,1)$ and setting $\ell=1$ yields the following
result.


\begin{corollary}\label{cp2}
Let $A=\{a_1,\ldots,a_k\}$ be any ordered subset of $\NN$.
Then the generating function $g_A(x;y)=\sum_{n\geq0}\sum_{\sigma\in
P_n^A}\levels(\sigma)x^ny^{\parts(\sigma)}$ is given by
$$\frac{y^2\left(1-y^2\sum_{i=1}^kx^{2a_i}\right)\sum_{i=1}^kx^{2a_i}(1+2x^{a_i}y)
+2y^4\left(1+y\sum_{i=1}^kx^{a_i}\right)\sum_{i=1}^kx^{4a_i}}
{\left(1-y^2\sum_{i=1}^kx^{2a_i}\right)^2}.$$
\end{corollary}

For example, applying Corollary~\ref{cp2} with $A=\NN$ gives that the generating
function $g_{\NN}(x;y)$ for the number of levels in all palindromic compositions of
$n$ with $m$ parts in $\NN$ is given by
\begin{equation} \label{levelsgf}
\frac{x^2y^2\biggl(2x^4(x+1)y^3+x^2(1-3x^2)(1+x+x^2)y^2+(1-x^4)(2x(1+x)y+1+x+x^2)\biggr)}
{(1+x^2)(1+x+x^2)(1-x^2-x^2y^2)^2}.
\end{equation}
Rewriting $\frac{1}{(1-x^2-x^2y^2)^2}$ as
$$\frac{1}{(1-x^2)^2(1-\frac{x^2y^2}{1-x^2})^2}=\frac{1}{(1-x^2)^2}\sum_{m\ge
0}(m+1)\frac{x^{2m}}{(1-x^2)^m}y^{2m}$$ allows us to compute the
generating function $l_m(x)$ for the number of levels in
palindromic compositions of $n$ with a given number of parts, $m$,
by looking at the coefficient of $y^m$ in
expression~(\ref{levelsgf}):
$$l_m(x)=\left\{\begin{array}{lll}
\frac{x^2}{1-x^2}&\mbox{ for }& m=2\\
\frac{(2m'-1-(2m'-3)x^2)x^{2m'}}{(1+x^2)(1-x^2)^{m'}}& \mbox{ for } & m=2m', \: m'\ge 2\\
\frac{2(1+x)(m'+(m'-1)x+m'x^2)x^{2m'+1}}{(1+x^2)(1+x+x^2)(1-x^2)^{m'}}& \mbox{ for } & m=2m'+1, \: m'\ge 1
\end{array}.\right.
$$
In addition, setting $y=1$  in (\ref{levelsgf}) gives that the
generating function for the number of levels in the palindromic compositions
of $n$ with parts in $\NN$ (see~\cite{CGH}, Theorem~6) is given by
    $$g_{\NN}(x;1)=\frac{x^2(1+3x+4x^2+x^3-x^4-4x^5-6x^6)}{(1+x^2)(1+x+x^2)(1-2x^2)^2}.$$
If we let $A=\{1,k\}$ in Corollary~\ref{cp2}, then we get that
$g_{\{1,k\}}(x;y)$ is given by
$$\frac{y^2(x^2+x^{2k})+2y^3(x^3+x^{3k})+y^4(x^4+x^{4k}-2x^{2(k+1)})+2y^5(x^{k+4}-x^{2k+3}-x^{3k+2}+x^{4k+1})}{(1-y^2x^2-y^2x^{2k})^2}.$$
Setting $y=1$ in the above expression yields that the generating
function for the number of levels in the palindromic compositions
of $n$ with any number of parts in $\{1,k\}$ is given by
$$g_{\{1,k\}}(x;1)=\frac{x^2+x^{2k}+x^3+x^{3k}+x^4+x^{4k}+2(x^{k+4}-x^{2(k+1)}-x^{2k+3}-x^{3k+2}+x^{4k+1})}{(1-y^2x^2-y^2x^{2k})^2}.$$
This result was not explicitly stated in~\cite{CH2}, but can be
easily computed from the generating functions for other quantities
given in~\cite{CH2}.

We look next at $\Aodd$. Applying Corollary~\ref{cp2} for this
case, we get that
$$g_{A}(x;y)=\frac{y^2\left(1-\frac{x^2y^2}{1-x^4}\right)\left(\frac{x^2}{1-x^4}+2\frac{x^3y}{1-x^6}\right)+2\frac{x^4y^4}{1-x^8}\left(1+\frac{x\,y}{1-x^2}\right)}{\left(1-\frac{x^2y^2}{1-x^4}\right)^2}.$$
Furthermore, if we set $y=1$ in the above expression, then we get
that the generating function for the number of levels in the
palindromic compositions of $n$ with any number of odd parts is
given by
$$g_A(x;1)=\frac{x^2(1+2x+2x^2+2x^3+2x^4+2x^5-2x^6+2x^7-4x^8-2x^9-4x^{10}-2x^{11}-x^{12})}{(1+x^4)(1-x^2-x^4)^2(1+x^2+x^4)},$$
which extends the work of Grimaldi~\cite{G}.

Finally, applying Corollary~\ref{cp2} for $A=\NN-\{k\}$ gives that the generating function
$g_{\NN-\{k\}}(x;y)$ is given by
$$\frac{y^2\left(1-\frac{y^2x^2}{1-x^2}+y^2x^{2k}\right)\left(\frac{x^2}{1-x^2}-x^{2k}+\frac{2yx^3}{1-x^3}-2yx^{3k}\right)
+2y^4\left(1+\frac{yx}{1-x}-yx^{k}\right)\left(\frac{x^4}{1-x^4}-x^{4k}\right)}
{\left(1-\frac{y^2x^2}{1-x^2}+y^2x^{2k}\right)^2}.$$ In
particular, when setting $y=1$ in the above expression we get that
the generating function for the number of
levels in the palindromic compositions of $n$ with any number of
parts in $A=\NN-\{k\}$ is given by
\begin{eqnarray*}
\frac{x^2(1+3x+4x^2+x^3-x^4-4x^5-6x^6)+x^{2k}(x^4-1)(1+x-2x^2-5x^3-5x^4)}{(1+x^2)(1+x+x^2)(1-2x^2+x^{2k}-x^{2(k+1)})^2}\\
+\frac{(x^2-1)(2x^{k+4}+2x^{3k}(1+x^2)(1-2x^2)+x^{4k}(1+x)(3-x)(1+x^2))}{(1+x^2)(1-2x^2+x^{2k}-x^{2(k+1)})^2}.
\end{eqnarray*}
This extends the work of Chinn and Heubach~\cite{CH}. Likewise, we
can extend the work of Grimaldi~\cite{G2} by setting $k=1$ to get that
$$g_{\NN-\{1\}}(x;1)=\frac{(1+x+3x^2+2x^3-5x^6-3x^7-x^8)x^4}{(1-x^2-x^4)^2(1+x^2)(1+x+x^2)}.$$

\subsection{Carlitz Compositions with parts in $A$}\label{seccar1} A
{\em Carlitz composition} of $n$, introduced in~\cite{C}, is a
composition of $n$ in which no adjacent parts are the same. In
other words, a Carlitz composition $\sigma$ is a composition with
$\levels(\sigma)=0$. We will derive results on the set of Carlitz
compositions of $n$ with parts in $A$, denoted by $E_n^A$. In
this section we study the generating functions for the number of
Carlitz compositions of $n$ with parts in $A$ with respect to the
number of rises and drops.

\subsubsection{Number of Carlitz compositions}
We denote the generating function for the number of Carlitz
compositions of $n$ with $m$ parts in $A$ with respect to the number of
rises and drops by $E_A(x;y;r,d)$, that is,
$$E_A(x;y;r,d)=\sum_{n\geq0}\sum_{\sigma\in E_n^A}x^ny^{\parts(\sigma)}r^{\rises(\sigma)}d^{\drops(\sigma)}.$$
Note that $E_A(x;y;r,d)=C_A(x;y;r,0,d)$. Therefore, Theorem~\ref{mth}(i) for $\ell=0$ gives the following
result.

\begin{corollary}\label{car1}
Let $A=\{a_1,\ldots,a_k\}$ be any ordered subset of $\NN$. Then
$$E_A(x;y;r,d)=1+\dfrac{\ds \sum_{j=1}^k\left(\dfrac{x^{a_j}y}{1+x^{a_j}y\,d}\ds \prod_{i=1}^{j-1}\dfrac{1+x^{a_i}y\,r}{1+x^{a_i}y\,d}\right)}
{1-\ds \sum_{j=1}^k\left(\dfrac{x^{a_j}y\,d}{1+x^{a_j}y\,d} \ds \prod_{i=1}^{j-1}\dfrac{1+x^{a_i}y\,r}{1+x^{a_i}y\,d}\right)}.$$
\end{corollary}

For example, with $r=d=1$, Corollary~\ref{car1} gives that
the generating function for the number of Carlitz compositions with
$m$ parts in $A$ (for the case $A=\NN$, see~\cite{C}) is given by
$$E_A(x;y;1,1)=\dfrac{1}{1-\sum_{j=1}^k\dfrac{x^{a_j}y}{1+x^{a_j}y}}.$$

Applying Corollary~\ref{car1} for $A=\{a,b\}$ and $r=d=1$ yields the generating function for the number of Carlitz compositions of $n$ with $m$ parts in $\{a,b\}$ is given by
$$\frac{(1+x^ay)(1+x^by)}{1-x^{a+b}y^2}=1+(x^a+x^b)y+\sum_{m\ge 1}x^{m(a+b)}(2y^{2m}+(x^a+x^b)y^{2m+1}).$$
In particular, setting $y=1$ in the expression above yields that the generating function for the number of Carlitz compositions of $n$ with parts in $\{a,b\}$
is given by
$$\frac{(1+x^a)(1+x^b)}{1-x^{a+b}}.$$

{\bf Remark:} In the case $A=\{a,b\}$, the requirement that no adjacent parts are to be the same restricts the compositions to those with alternating $a$'s and $b$'s. This results in the following possibilities:
\begin{eqnarray} \label{CC} n \qquad \quad &\mbox{Carlitz compositions of }n\nonumber \\
n'(a+b)\quad &abab \ldots ab \mbox{ and }baba \dots ba \nonumber \\
n'(a+b)+a & abab \ldots aba\\ \nonumber
n'(a+b)+b & babab \ldots ab
\end{eqnarray}
Thus, the number of Carlitz compositions of $n>0$ is 2 if
$n=n'(a+b)$, 1 if $n=n'(a+b)+a$ or $n=n'(a+b)+b$, and 0 otherwise.

\subsubsection{Number of Rises and Drops}
We now study the number of rises (drops) in
all Carlitz compositions of $n$ with $m$ parts in $A$. Once more, the number of rises equals the number of drops. Using Corollary~\ref{car1} to find an
explicit expression for $\left.\frac{\partial }{\partial r}
E_A(x;y;r,1)\right|_{r=1}$ gives
the following result.

\begin{corollary}\label{car2}
Let $A=\{a_1,\ldots,a_k\}$ be any ordered subset of $\NN$. Then
the generating functions $\sum_{n\geq0}\sum_{\sigma\in
E_n^A}\rises(\sigma)x^ny^{\parts(\sigma)}$ and $\sum_{n\geq0}\sum_{\sigma\in
E_n^A}\drops(\sigma)x^ny^{\parts(\sigma)}$ are given by
$$\dfrac{\ds \sum_{j=1}^k\left(\dfrac{x^{a_j}y}{1+x^{a_j}y}\ds \sum_{i=1}^{j-1}\dfrac{x^{a_i}y}{1+x^{a_i}y}\right)}
{\left(1-\ds \sum_{j=1}^k\dfrac{x^{a_j}y}{1+x^{a_j}y}\right)^2}.$$
\end{corollary}

Setting $A=\NN$ and $y=1$ in Corollary~\ref{car2} yields that the generating function for the number of rises (drops) in the Carlitz compositions of $n$ with parts in $\NN$ is given
by
$$\frac{\sum_{j\geq1}\left(\frac{x^j}{1+x^j}\sum_{i=1}^{j-1}\frac{x^i}{1+x^i}\right)}
{\left(1-\sum_{j\ge1}\frac{x^j}{1+x^j}\right)^2}.$$

Applying Corollary~\ref{car2} for $A=\{a,b\}$ gives that
\begin{eqnarray*}
\sum_{n\geq0}\sum_{\sigma\in
E_n^A}\rises(\sigma)x^ny^{\parts(\sigma)}&=&
\frac{x^{a+b}y^2(1+x^ay)(1+x^by)}{(1-x^{a+b}y^2)^2}\\
&=&\sum_{m\geq1}
x^{m(a+b)}\left((2m-1)y^{2m}+m(x^a+x^b)y^{2m+1}\right),
\end{eqnarray*}
where the second equation follows after collecting even and odd powers of $y$.

In particular, setting $y=1$ in the expression above yields that the generating function for the number of rises (drops) in the Carlitz
compositions of $n$ with parts in $\{a,b\}$ is given by
$$\frac{x^{a+b}(1+x^a)(1+x^b)}{(1-x^{a+b})^2}.$$

Thus, the number of rises (drops) in Carlitz compositions of $n\ge (a+b)$ with parts in $\{a,b\}$ is given by
$$ n' \:\mbox{  if  }\: n=(a+b)n'+a \mbox{  or  } n=(a+b)n'+b\quad \mbox{and} \quad 2n'-1 \:\mbox{  if  }\: n=(a+b)n' \quad\mbox{for } n' \ge 1.$$
This follows immediately from (\ref{CC}) since there is a rise for every occurrence of ``$ab$''. If $n=(a+b)n'$ and the composition starts with $a$, then there are $n'$ rises. For the composition that starts with $b$, there is one less rise, for a total of $2n'-1$ rises. If $n$ is not a multiple of $a+b$, then the composition starts with $r$, where $n=(a+b)n'+r$. In either case, there are exactly $n'$ rises, as there are $n'$ occurrences of ``$ab$'' in the composition.


\subsection{Carlitz palindromic compositions}\label{seccarp}
A {\em Carlitz palindromic composition} of $n$ is both a Carlitz composition and a palindromic
composition. Let $F_n^A=E_n^A\cap
P_n^A$ be the set of all Carlitz palindromic compositions of $n$
with parts in $A$.

\subsubsection{ Number of Carlitz palindromic compositions}
We denote the generating function for the number of Carlitz
palindromic compositions of $n$ with $m$ parts in $A$ with respect
to the number of rises by $F_A(x;y;r)$, that is,
$$F_A(x;y;r)=\sum_{n\geq0}\sum_{\sigma\in F_n^A}x^ny^{\parts(\sigma)}r^{\rises(\sigma)}.$$
Note that $F_A(x;y;r)=P_A(x;y;r,0,1)$. Using Theorem~\ref{mth}(ii) for $\ell=0$ and $d=1$ gives the following result.

\begin{corollary}\label{carp1}
Let $A=\{a_1,\ldots,a_k\}$ be any ordered subset of $\NN$. Then
$$F_A(x;y;r)=1+\dfrac{\sum\limits_{i=1}^k\dfrac{x^{a_i}y}{1+x^{2a_i}y^2r}}
{1-\sum\limits_{i=1}^k\dfrac{x^{2a_i}y^2r}{1+x^{2a_i}y^2r}}.$$
\end{corollary}

Applying Corollary~\ref{carp1} for $A=\{a,b\}$ and $y=r=1$ yields that the generating function for the number of Carlitz palindromic compositions of $n$ with parts in $\{a,b\}$ is given by
$$\frac{1+x^a+x^b-x^{a+b}}{1-x^{a+b}}.$$
Thus, the number of Carlitz palindromic compositions of $n$ with parts in $\{a,b\}$ is 1 if $n =(a+b)n'+a$ or $n =(a+b)n'+b$ for some $n'\ge 0$, and 0 otherwise. This follows immediately from (\ref{CC}), since the Carlitz compositions for $n=(a+b)n'$ are not symmetric.

\subsubsection{Number of Rises and Drops}
We now study the number of rises (drops) in
all Carlitz palindromic compositions of $n$ with $m$ parts in $A$. Using Corollary~\ref{carp1} to compute $\left.\frac{\partial }{\partial r} F_A(x;y;r)\right|_{r=1}$ gives the following result.

\begin{corollary}\label{carp2}
Let $A=\{a_1,\ldots,a_k\}$ be any ordered subset of $\NN$.
Then the generating function for the number of rises in
all Carlitz palindromic compositions of $n$ with $m$ parts in $A$
is given by
$$\left.\frac{\partial }{\partial r} F_A(x;y;r)\right|_{r=1}
=\dfrac{\sum_{i=1}^k\frac{x^{3a_i}y^3}{(1+x^{2a_i}y^2)^2}\left(\sum_{i=1}^k\frac{x^{2a_i}y^2}{1+x^{2a_i}y^2}-1\right)+
\sum_{i=1}^k\frac{x^{a_i}y}{1+x^{2a_i}y^2}\sum_{i=1}^k\frac{x^{2a_i}y^2}{(1+x^{2a_i}y^2)^2}}
{\left(1-\sum\limits_{i=1}^k\frac{x^{2a_i}y^2}{1+x^{2a_i}y^2}\right)^2}.$$
\end{corollary}

Applying Corollary~\ref{carp2} for $A=\{a,b\}$ gives that
$$\sum_{n\geq0}\sum_{\sigma\in
F_n^A}\rises(\sigma)x^ny^{\parts(\sigma)}=
\frac{x^{a+b}y^3(x^a+x^b)}{(1-x^{a+b}y^2)^2}
=(x^a+x^b)\sum_{m\geq1}
m \, x^{m(a+b)}y^{2m+1}.$$

In particular, setting $y=1$ in the expression above yields that the generating function for the number of rises (drops) in all Carlitz
palindromic compositions of $n$ with parts in $\{a,b\}$ is given by
$$\frac{x^{a+b}(x^a+x^b)}{(1-x^{a+b})^2}.$$

Thus, the number of rises (drops) in the Carlitz palindromic compositions of $n \ge a+b$ with parts in $\{a,b\}$ is given by
$$ n' \: \mbox{  if  }\: n=(a+b)n'+a \mbox{  or  } n=(a+b)n'+b \: \mbox{  for  }\: n'\ge 1\quad \mbox{and}\quad  0 \:\mbox{  otherwise} .$$
This follows immediately from (\ref{CC}), as the Carlitz compositions for $n=(a+b)n'+a$ and $n=(a+b)n'+b$ are symmetric.


\subsection{Partitions with parts in $A$}\label{secpar}
A {\em partition} $\sigma$ of $n$ is a composition of $n$ with
$\rises(\sigma)=0$. Let $G_n^A$ be the set of all partitions of $n$
with parts in $A$.

\subsubsection{Number of partitions}
We denote the generating function for the number of partitions of $n$ with $m$ parts in $A$ with respect
to the number of levels and drops by
$$G_A(x;y;\ell,d)=\sum_{n\geq0}\sum_{\sigma\in
G_n^A}x^ny^{\parts(\sigma)}\ell^{\levels{\sigma}}d^{\drops(\sigma)}.$$
Note that $G_A(x;y;\ell,d)=C_A(x;y;0,l,d)$. Using Theorem~\ref{mth}(i) for $r=0$ we get the
following result.

\begin{corollary}\label{par1}
Let $A=\{a_1,\ldots,a_k\}$ be any ordered subset of $\NN$. Then
the generating function $G_A(x;y;\ell,d)$ is given by
$$1+\dfrac{\sum_{j=1}^k\left(\frac{x^{a_j}y}{1-x^{a_j}y(\ell-d)}\prod_{i=1}^{j-1}\frac{1-x^{a_i}y\ell}{1-x^{a_i}y(\ell-d)}\right)}
{1-d\sum_{j=1}^k\left(\frac{x^{a_j}y}{1-x^{a_j}y(\ell-d)}\prod_{i=1}^{j-1}\frac{1-x^{a_i}y\ell}{1-x^{a_i}y(\ell-d)}\right)}.$$
\end{corollary}

For example, if we apply Corollary~\ref{par1} for $A=\NN$ and $\ell=d=1$ and use the identity
\begin{equation}\label{eqd1}
\sum_{j=1}^k x^{a_j}\prod_{i=1}^{j-1}(1-x^{a_i}\alpha)=\frac{1}{\alpha}\left(1-\prod_{j=1}^k (1-x^{a_j}\alpha)\right),
\end{equation}
then we get that the generating function for the number of partitions of
$n$ with $m$ parts in $A=\NN$ is given by
$$F_\NN(x;y;1,1)=\prod_{j\geq1}(1-x^jy)^{-1}.$$

Note that the identity in (\ref{eqd1}) follows from the fact that
$$1-\alpha \sum_{j=1}^k x^{a_j}\prod_{i=1}^{j-1}(1-x^{a_i}\alpha)=\left(1-\prod_{j=1}^k (1-x^{a_j}\alpha)\right),$$
which can be easily proved by induction.

If we apply Corollary~\ref{par1} to $A=\{a,b\}$ and set $y=\ell=d=1$, then we get that the generating function for
the number of partitions of $n$ with parts in $A$ is given by
$$\frac{1}{1-x^a-x^b(1-x^a)}=\frac{1}{(1-x^a)(1-x^b)}.$$
In particular, if $A=\{1,k\}$ then we have that the number of
partitions of $n$ with parts in $A$ is given by $\lfloor (n+k)/k \rfloor$. This can be easily explained by the following observation. For $n\in [n'k,(n'+1)k)$, the only partitions are those consisting of all 1's, one $k$ and all 1's,\ldots,$n'$ $k$'s and all 1's, for a total of $n'+1=\lfloor (n+k)/k \rfloor$ partitions.

Another interesting example, namely setting $\ell=0$ and $d=1$ in
Corollary~\ref{par1}, gives that the generating function for the
number of partitions of $n$ with $m$ parts in $A$ in which no
adjacent parts are the same is given by
$$G_A(x;y;0,1)=\dfrac{1}{1-\sum_{j=1}^kx^{a_j}y\prod_{i=1}^j(1+x^{a_i}y)^{-1}}=\prod_{j=1}^k(1+x^{a_j}y), $$
where the second equality is easily proved by induction.
In particular, the generating function for the number of partitions
of $n$ with parts in $\NN$ in which no adjacent parts are the same
is given by $\prod_{j\geq1}(1+x^j)$.

\subsubsection{Number of levels and drops}
We now study the number of levels and drops in all partitions of $n$. Using Corollary~\ref{par1} to compute
$\left.\frac{\partial}{\partial
\ell}G_A(x;y;\ell,1)\right|_{\ell=1}$ and
$\left.\frac{\partial}{\partial d}G_A(x;y;1,d)\right|_{d=1}$,
 we get the following
result.

\begin{corollary}\label{partld}
Let $A=\{a_1,\ldots,a_k\}$ be any ordered subset of $\NN$. Then
the generating function $\sum_{n\geq0}\sum_{\sigma\in
G_n^A}\levels(\sigma)x^ny^{\parts(\sigma)}$ is given by
$$\frac{\sum_{j=1}^k\left(x^{2a_j}y^2\prod_{i=1}^{j-1}(1-x^{a_i}y)\right)
-\sum_{j=1}^k\left(x^{a_j}y\prod_{i=1}^{j-1}(1-x^{a_i}y)\sum_{i=1}^{j-1}\frac{x^{2a_i}y^2}{1-x^{a_i}y}\right)}
{\prod_{j=1}^k(1-x^{a_j}y)^2},$$ and the generating function
$\sum_{n\geq0}\sum_{\sigma\in
G_n^A}\drops(\sigma)x^ny^{\parts(\sigma)}$ is given by
$$\frac{\left(1-\prod_{j=1}^k(1-x^{a_j}y)\right)^2
-y^2\sum_{j=1}^k\left(x^{a_j}\prod_{i=1}^{j-1}(1-x^{a_i}y)\sum_{i=1}^{j}x^{a_i}\right)}
{\prod_{j=1}^k(1-x^{a_j}y)^2}.$$
\end{corollary}

\begin{proof}
We give a sketch of the proof for the first generating function. Since $G_A(x;y;l,1)=1+\dfrac{S(\ell)}{1-S(\ell)}$, where
$$S(\ell)=\sum_{j=1}^k\left(\frac{x^{a_j}y}{1-x^{a_j}y(\ell-1)}\prod_{i=1}^{j-1}\frac{1-x^{a_i}y\ell}{1-x^{a_i}y(\ell-1)}\right)=\sum_{j=1}^k\left(g_j(\ell)\prod_{i=1}^{j-1}f_i(\ell)\right),$$
we get that $$\frac{\partial}{\partial
\ell}G_A(x;y;\ell,1)=\dfrac{\frac{\partial}{\partial\ell}S(\ell)}{(1-S(\ell))^2}=\dfrac{\sum_{j=1}^k\frac{\partial}{\partial\ell}g_j(\ell)\prod_{i=1}^{j-1}f_i(\ell)+g_j(\ell)\frac{\partial}{\partial\ell}\prod_{i=1}^{j-1}f_i(\ell)}{(1-S(\ell))^2}.$$
Using Equation~(\ref{eqdiff}) gives that $\frac{\partial}{\partial\ell}\prod_{i=1}^{j-1}f_i(\ell)=-\prod_{i=1}^k f_i(\ell)\sum_{j=1}^k\frac{x^{2a_j}y^2}{(1-x^{a_j}y\ell)(1-x^{a_j}y(\ell-1))}$. Computing $\frac{\partial}{\partial\ell}g_j(\ell)$, setting $\ell = 1$ in the expression for $\frac{\partial}{\partial
\ell}G_A(x;y;\ell,1)$, then using Equation~(\ref{eqd1}) to simplify the denominator gives the stated result.
\end{proof}

Applying Corollary~\ref{partld} to $A=\{a,b\}$ gives that the generating function for the number of levels in the partitions of $n$ with $m$ parts in $\{a,b\}$ is given by
$$\frac{x^{2a}y^2(1-x^by)+x^{2b}y^2(1-x^ay)}{(1-x^ay)^2(1-x^by)^2}.$$
In particular, the generating function for the number of levels in the partitions of $n$ with parts in $\{1,2\}$ is given by
$$\frac{x^2(1-x^3)}{(1-x)^4(1+x)^2}.$$

>From the second part of  Corollary~\ref{partld} we get for
$A=\{a,b\}$ that the generating function for the number of drops
in the partitions of $n$ with $m$ parts in $\{a,b\}$ is given by
$$\frac{x^{a+b}y^2}{(1-x^ay)(1-x^by)}.$$
In particular, setting $y=1$ in the above expression yields that the generating function for the number of drops in all partitions of $n$ with parts in $\{1,k\}$ is given by
$$\frac{x^{k+1}}{(1-x)(1-x^k)}.$$
Thus, the number of drops in the partitions of $n$ with parts in $\{1,k\}$ is $\lfloor (n-1)/k \rfloor$. This again follows from the specific structure of the partitions with parts in $\{1,k\}$. A single drop occurs in all the partitions that do not consist of either all 1's or all $k$'s. Thus, for $n \in [n'k+1,(n'+1)k)$, there are exactly $n'=\lfloor (n-1)/k \rfloor$ drops.

\section{Concluding Remarks}
We have provided a very general framework for answering questions
concerning the number of compositions, number of parts, and number
of rises, levels and drops in all compositions of $n$ with parts in
$A$. We have used this framework to investigate compositions, palindromic
compositions, Carlitz compositions, Carlitz palindromic
compositions and partitions of $n$. Our results generalize work by
several authors, and we have applied our results to the specific
sets studied previously, which has led to several new results. In
addition, our results can be applied to any set $A \subseteq \NN$,
which will allow for further study of special cases.

In addition, the techniques used in this paper can be used to
investigate products among the number of rises, levels and drops
which show interesting connections to the Fibonacci sequence, one
of the reasons Alladi and Hoggatt investigated the these
quantities for compositions with summands 1 and 2. For example, by
computing the derivative with respect to $d$ twice in Theorem 1.1
(ii) and setting $y=r=\ell=1$, we get that
$$\sum_{n\geq0}\sum_{\sigma\in C_n^A}\drops(\sigma)(\drops(\sigma)-1)x^n=\frac{2x^6}{(1-x-x^2)^3}=2x^3\sum_{n\geq3}\left(\sum_{a+b+c=n}F_aF_bF_c\right)x^n,$$
i.e., a convolution of three Fibonacci sequences. However, the
formulas for the various products become more complicated, and not
as easy to evaluate.

\end{document}